\documentclass{article}
\usepackage[affil-it]{authblk}
\usepackage[utf8]{inputenc}
\usepackage{amsmath, amsfonts, amsthm, amssymb, latexsym}
\usepackage{geometry}
\usepackage{enumerate}
\usepackage{multirow}
\usepackage{rotating}
\usepackage{graphicx,color}
\usepackage{float}
\usepackage{mathrsfs}

\newtheorem{thm}{Theorem}[section]

\newtheorem{define}[thm]{Definition}
\newtheorem{prop}[thm]{Proposition}

\title{An inverse problem for the fractional porous medium equation}
\author{Li Li}
\affil{Institute for Pure and Applied Mathematics, University of California,\\
Los Angeles, CA 90095, USA}
\date{}

\begin{document}
\maketitle

\noindent \textbf{ABSTRACT.}\, We consider a time-independent variable coefficients fractional porous medium equation and formulate an associated inverse problem. We determine both the conductivity and the absorption coefficient from exterior partial measurements of the Dirichlet-to-Neumann map. Our approach relies on a 
time-integral transform technique as well as the unique continuation property of the fractional operator.

\section{Introduction}
The classical porous medium equation
\begin{equation}\label{CPME}
\partial_t u- \Delta(|u|^{m-1}u)= 0,\qquad m>1
\end{equation}
appears in models for gas flow through porous media,
high-energy physics, population dynamics and many other contexts.
As a paradigm for nonlinear degenerate diffusion equations, the classical porous medium equation has been extensively studied so far. See \cite{vazquez2007porous}
for a survey in this field.

Results for Inverse problems related with the classical porous medium equation have been obtained in \cite{carstea2021inverse, carsteaU2021inverse}. There the authors studied time-independent variable coefficients porous medium equations and defined the associated Dirichlet-to-Neumann maps. It has been shown in \cite{carstea2021inverse, carsteaU2021inverse} that the coefficients in the equations can be uniquely determined from the knowledge of Dirichlet-to-Neumann maps. See \cite{carstea2021calder, kian2020recovery, krupchyk2019partial, shankar2020recovering} for other recent works on inverse problems for nonlinear elliptic and parabolic equations.

In recent years, nonlocal operators including the fractional Laplacian have attracted much attention. Replacing the classical Laplacian $- \Delta$ by the fractional Laplacian $(- \Delta)^s$ is motivated by the
need to describe processes involving anomalous diffusion. Such processes have been widely observed, making fractional Laplacian broadly applicable in fluid mechanics, biology, finance and many other disciplines.

In this paper, we will study an inverse problem related with
the following basic fractional porous medium equation
\begin{equation}\label{FPME}
\partial_t u+ (- \Delta)^s(|u|^{m-1}u)= 0,\qquad m>1,\quad 0<s<1,
\end{equation}
which is a natural combination of fractional diffusion and porous medium nonlinearities. Roughly speaking, this model describes anomalous diffusion through porous media.
See \cite{vazquez2014recent} and the references there for more background information.

To formulate our inverse problem, we first need to study an initial exterior problem for a time-independent variable coefficients fractional porous medium equation involving an absorption term.

More precisely, we consider the initial exterior problem 
\begin{equation}\label{fracPMie}
\left\{
\begin{aligned}
\partial_t u+ L^s_{\gamma}(|u|^{m-1}u)+ \lambda(x)u &= 0,\quad \Omega\times (0, T),\\
u&= g,\quad \Omega_e\times (0, T),\\
u&= 0,\quad \Omega\times \{0\}\\
\end{aligned}
\right.
\end{equation}
where $\Omega_e:= \mathbb{R}^n\setminus \bar{\Omega}$, $g$ satisfies $\mathrm{supp}\,g\subset
\Omega_e\times [0, T]$ and the elliptic operator $L_{\gamma}$ is defined by
\begin{equation}\label{Lgamma}
L_{\gamma}:= -\mathrm{div}(\gamma(x)\nabla).
\end{equation}
Here we assume $\lambda\in C^\infty(\bar{\Omega})$,
$0< \gamma\in C^\infty(\mathbb{R}^n)$ and $\gamma= 1$ in $\Omega_e$ for convenience. Note that the fractional operator $L^s_{\gamma}$ coincides with
$(- \Delta)^s$ when $\gamma= 1$ in $\mathbb{R}^n$.

Formally we define the Dirichlet-to-Neumann map $\Lambda_{\gamma,\lambda}$ associated with (\ref{fracPMie}) by
\begin{equation}\label{glDN}
\Lambda_{\gamma,\lambda}: g\to L^s_{\gamma}(|u|^{m-1}u)|_{\Omega_e\times (0, T)}.
\end{equation}
We will show that (\ref{glDN}) is well-defined at least for
time-independent exterior data $g:= g_0(x)\in S_W$. Here the set  $S_W$ is defined by
\begin{equation}\label{setSW}
S_W:=\{g_0:\,|g_0|^{m-1}g_0\in C^\infty_c(W)\}
\end{equation}
where $W\subset \Omega_e$ is nonempty and open.

Our goal is to determine the absorption coefficient $\lambda$ and the conductivity $\gamma$ in $\Omega$ from partial measurements of $\Lambda_{\gamma,\lambda}$.
The following theorem is the main result in this paper.

\begin{thm}
Let $W_1, W_2\subset \Omega_e$ be nonempty and open.
Suppose $\lambda^{(1)}, \lambda^{(2)}\in C^\infty(\bar{\Omega})$,
$0< \gamma^{(1)}, \gamma^{(2)}\in C^\infty(\mathbb{R}^n)$ and $\gamma^{(1)}= \gamma^{(2)}= 1$ in $\Omega_e$.
Suppose 
\begin{equation}\label{IdDN}
\Lambda_{\gamma^{(1)},\lambda^{(1)}}g|_{W_2\times (0, T)}
= \Lambda_{\gamma^{(2)},\lambda^{(2)}}g|_{W_2\times (0, T)}
\end{equation}
for all $g\in S_{W_1}$. Then $\gamma^{(1)}= \gamma^{(2)}$
and $\lambda^{(1)}= \lambda^{(2)}$ in $\Omega$.
\end{thm}

We remark that our problem can be viewed as a fractional analogue of the inverse problems studied in \cite{carstea2021inverse, carsteaU2021inverse}. Our problem can also be regarded as a nonlinear parabolic variant of the inverse problem studied in 
\cite{ghosh2021calder} where the authors considered the
exterior Dirichlet problem
\begin{equation}\label{eDp}
L^s_{\gamma}u= 0\,\,\, \text{in}\,\,\Omega,\qquad
u= g\,\,\,\, \text{in}\,\,\Omega_e.
\end{equation}
It has been proved in \cite{ghosh2021calder} that $\gamma$ in $\Omega$ can be determined from exterior partial measurements of the associated Dirichlet-to-Neumann map 
\begin{equation}\label{DNlin}
\Lambda^{lin}_\gamma: g\to L^s_{\gamma}u|_{\Omega_e}.
\end{equation}
Later we will see how to relate our nonlinear fractional parabolic problem to this linear fractional elliptic problem via a time-integral transform.

We also mention that the study of inverse problems for space-fractional operators was initiated in \cite{ghosh2020calderon} and further results in this direction can be found in many recent works.
Stability and single measurement results for the fractional Calder\'on problem have been obtained in \cite{ghosh2020uniqueness, ruland2020fractional}.
See \cite{bhattacharyya2021inverse, covi2021uniqueness, covi2020higher} for inverse problems for
fractional Schr\"odinger operators with local, quasilocal and nonlocal perturbations. See \cite{covi2020inverse, lai2021inverse, li2021determining} for inverse problems for
fractional operators in the magnetic setting. See \cite{li2021inverse} for inverse problems for fractional operators in the elastic setting.
See \cite{covi2021calder} for inverse problems for directionally antilocal operators. See \cite{lai2020calderon, li2021fractional} for inverse problems for linear fractional parabolic operators. See \cite{li2021inversediff} for an inverse problem for a different nonlinear fractional parabolic operator.

All the uniqueness results obtained in the works mentioned above rely on exploiting unique continuation properties of associated fractional operators. This typical nonlocal phenomenon makes inverse problems for fractional operators often more manageable than their classical counterparts. Later we will see how to determine $\lambda$ in an immediate way by using the unique continuation property of $L^s_{\gamma}$ after we determine $\gamma$.

The rest of this paper is organized in the following way. We summarize the background knowledge in Section 2.
We study the initial exterior problem (\ref{fracPMie}) in Section 3 via the theory of monotone operators in Hilbert spaces. In Section 4,
we first use a time-integral transform technique to determine 
the conductivity $\gamma$; Then we use the unique continuation property to determine the absorption coefficient $\lambda$.
\medskip

\noindent \textbf{Acknowledgements.} The author would like to thank Professor Gunther Uhlmann for suggesting the problem and for helpful discussions.

\section{Preliminaries}
Throughout this paper we use the following notations.
\begin{itemize}
\item Fix the space dimension $n\geq 2$

\item Fix the fractional powers $0< s< 1$ and $m> 1$.

\item For convenience, we write $u^m:= |u|^{m-1}u$;
$u:= v^{\frac{1}{m}}$ if $v= u^m$.

\item Fix the constant $T> 0$ and $t$ denotes the time variable.

\item $\Omega$ denotes a bounded domain with smooth boundary and
$\Omega_e:= \mathbb{R}^n\setminus\bar{\Omega}$.

\item Suppose $u$ is an $(n+1)$-variable function. Then 
$u(t)$ denotes the $n$-variable function $u(\cdot, t)$.

\item $c, C, C', C_1,\cdots$ denote positive constants. We write $C_I$ to emphasize the dependence on the parameter $I$.

\item $X^*$ denotes the dual space of X.
$\langle\cdot, \cdot\rangle$ denotes the dual pairing.
\end{itemize}

\subsection{Function spaces}
Throughout this paper we refer all function spaces to real-valued
function spaces.

We use $H^r$ to denote $W^{r,2}$-type Sobolev spaces.

Let $U$ be an open set in $\mathbb{R}^n$. Let $F$ be a closed set in $\mathbb{R}^n$. Then
$$H^r(U):= \{u|_U: u\in H^r(\mathbb{R}^n)\},\qquad 
H^r_F(\mathbb{R}^n):= 
\{u\in H^r(\mathbb{R}^n): \mathrm{supp}\,u\subset F\},$$
$$\tilde{H}^r(U):= 
\mathrm{the\,\,closure\,\,of}\,\, C^\infty_c(U)\,\,\mathrm{in}\,\, H^r(\mathbb{R}^n).$$

For $r\in\mathbb{R}$, we have the natural identifications
$$H^{-r}(\mathbb{R}^n)= H^r(\mathbb{R}^n)^*,\quad \tilde{H}^r(\Omega)= H^r_{\bar{\Omega}}(\mathbb{R}^n),\quad 
H^{-r}(\Omega)= \tilde{H}^r(\Omega)^*.$$

Let $X$ be a Banach space. We use $C([0, T]; X)$
to denote the space consisting of the corresponding Banach space-valued continuous functions on $[0, T]$. 
$L^2(0, T; X)$ (resp., $H^1(0, T; X)$) denotes the space consisting of the corresponding Banach space-valued $L^2$-functions (resp., $H^1$-functions).

\subsection{Fractional operators}
We briefly present the precise definition of $L^s_{\gamma}$
via the semigroup approach. We also list some of its most important properties for later use.

It is well-known that there exists a unique symmetric heat kernel $p_t$ s.t. $p_t(x, y)$ is smooth jointly in
$t>0, x, y\in \mathbb{R}^n$ and
$$(e^{-t L_{\gamma}}f)(x)= \int p_t(x,y)f(y)\,dy,\qquad x\in \mathbb{R}^n,\, t>0,\,\,
f\in L^2(\mathbb{R}^n).$$
Moreover, we have the following Gaussian bounds
$$c_1e^{-c'_1\frac{|x-y|^2}{t}}t^{-\frac{n}{2}}\leq p_t(x, y)
\leq c_2e^{-c'_2\frac{|x-y|^2}{t}}t^{-\frac{n}{2}},\qquad
x, y\in \mathbb{R}^n,\,\,t>0.$$
We define
$$K(x, y):= C\int^\infty_0 p_t(x, y)\frac{dt}{t^{1+s}}.$$
Then we get the estimate
$$ \frac{C_1}{|x-y|^{n+2s}}\leq K(x, y)= K(y, x)\leq \frac{C_2}{|x-y|^{n+2s}},\qquad\quad x,y\in \mathbb{R}^n.$$
The fractional operator $L^s_{\gamma}$ is defined by
$$L^s_{\gamma}:= \frac{1}{\Gamma(-s)}\int^\infty_0(e^{-t L_{\gamma}}-\mathrm{Id})
\frac{dt}{t^{1+s}}$$
where $\Gamma$ is the Gamma function. It has been shown that
$$\langle L^s_{\gamma}u, v\rangle = \frac{1}{2}\iint
(u(x)- u(y))(v(x)- v(y))K(x, y)\,dxdy,\quad u,v\in H^s(\mathbb{R}^n).$$

It follows from Poincar\'e's inequality for the fractional Laplacian and the Lax-Milgram Theorem that the map
$$u\to L^s_{\gamma}u|_{\Omega}$$
gives a homeomorphism from $\tilde{H}^{s}(\Omega)$ to $H^{-s}(\Omega)$. Moreover, we have the following well-posedness
result.
\begin{prop}
For $f\in H^{-s}(\Omega)$ and $g\in H^s(\mathbb{R}^n)$, the 
exterior Dirichlet problem
\begin{equation}\label{eDpfg}
L^s_{\gamma}u= f\,\,\,\,\, \text{in}\,\,\,\Omega,\qquad
u= g\,\,\,\,\,\,\text{in}\,\,\,\Omega_e.
\end{equation}
has a unique solution $u\in H^s(\mathbb{R}^n)$ and 
$$||u||_{H^s(\mathbb{R}^n)}\leq C(||f||_{H^{-s}(\Omega)}+
||g||_{H^s(\mathbb{R}^n)}).$$
\end{prop}

The following proposition is the unique continuation property of 
$L_{\gamma}$. Its proof is based on the Caffarelli-Silvestre definition of
the fractional Laplacian introduced in \cite{caffarelli2007extension}.
\begin{prop}\label{fracUCP}
Let $u\in H^{s}(\mathbb{R}^n)$. Let $W$ be open. Suppose 
$$L_{\gamma}u= u= 0\quad\text{in}\,\,W.$$
Then $u= 0$ in $\mathbb{R}^n$.
\end{prop}

In fact, all the results in this subsection hold true for the general variable coefficients operator
$$L_A:= -\mathrm{div}(A(x)\nabla)$$
where the smooth real symmetric matrix-valued function $A(x):= (a_{i,j}(x))$ satisfies the uniformly elliptic condition, i.e.
$$C^{-1}_A|\xi|^2\leq \sum_{1\leq i,j\leq n}a_{i, j}(x)\xi_i\xi_j\leq C_A|\xi|^2,\qquad x, \xi\in \mathbb{R}^n.$$
We refer readers to \cite{ghosh2021calder} or \cite{ghosh2017calderon} for more details.

\section{Forward problem}
To study (\ref{fracPMie}), we make the substitution $w:= u-g$ so
$w= 0$ in $\Omega_e$. 
Note that $u^m= w^m+ g^m$ since the supports of $w, g$ are disjoint so (\ref{fracPMie}) can be converted into the initial value problem
\begin{equation}\label{fracPMivp}
\left\{
\begin{aligned}
\partial_t w+ L^s_{\gamma}(w^m)+ \lambda(x)w &= f,\quad \Omega\times (0, T),\\
w&= 0,\quad \Omega\times \{0\}.\\
\end{aligned}
\right.
\end{equation}

Let us briefly present some basic concepts
in the theory of monotone operators in Hilbert spaces.

\begin{define}
Let $X$ be a Banach space. Let $\Psi$ be a convex, lower semicontinuous functional on $X$, we say that a multivalued map $A$ from $X$ to $X^*$ is the 
subdifferential of $\Psi$ if for each $z\in X$,
$$Az= \{z^*\in X^*: \Psi(z)-\Psi(z')\leq \langle z^*, z-z'\rangle, \forall \,z'\in X\}.$$
\end{define}

Here we are only interested in single-valued maps defined in Hilbert spaces. 

It is well-known that the subdifferential of a lower semicontinuous convex functional is maximal monotone. Recall that for an operator $A$ defined in a Hilbert space $H$, $A$ is maximal monotone if and only if $A$ is monotone and the range of
Id$+A$ is $H$. See \cite{barbu2010nonlinear} for more details.

In our case, we consider the functional $\Psi$ defined on $H^{-s}(\Omega)$ given by
$$\Psi(u):= \frac{1}{m+ 1}\int_{\Omega}|u|^{m+1}$$
for $u\in H^{-s}(\Omega)\cap L^1(\Omega)$ satisfying $u^{m+1}\in L^1(\Omega)$
and $\Psi(u):= +\infty$ otherwise. 

We define the operator $A$ in $H^{-s}(\Omega)$ by
$$Az:= L^s_{\gamma}(z^m)|_{\Omega}$$
with the domain
$$D(A):= \{z\in H^{-s}(\Omega)\cap L^1(\Omega): z^m\in \tilde{H}^s(\Omega)\}.$$

The following result is Proposition 3.1 in \cite{bonforte2014existence}.
\begin{prop}
$\Psi$ is a convex, lower semicontinuous functional on $H^{-s}(\Omega)$; $A$
is the subdifferential of $\Psi$ and thus a maximal monotone operator in 
$H^{-s}(\Omega)$.
\end{prop}

Now by Theorem 4.11 and Remark 4.5 in \cite{barbu2010nonlinear}, we get the following existence, uniqueness and regularity results for (\ref{fracPMivp}).
\begin{prop}
Let $f\in L^2(0, T; H^{-s}(\Omega))$. Then there exists a unique solution
$$w\in C([0, T]; H^{-s}(\Omega))\cap H^1(0, T; H^{-s}(\Omega))$$
of (\ref{fracPMivp}). Moreover, $w(t)\in D(A)$ for $t\in (0, T)$.
\end{prop}

Hence, at least for $g$ satisfying $\mathrm{supp}\,g\subset 
\Omega_e\times [0, T]$ and 
$$f:= -L^s_{\gamma}(g^m)|_{\Omega\times (0, T)}\in L^2(0, T; H^{-s}(\Omega)),$$
$u:= w_f+ g$ gives the unique solution of (\ref{fracPMie}) and the Dirichlet-to-Neumann map $\Lambda_{\gamma,\lambda}$ given by (\ref{glDN}) is well-defined for such $g$.

\section{Inverse problem}
To study the inverse problem, we will focus on the time-independent exterior data $g\in S_W$ where the set $S_W$ is defined by (\ref{setSW}).

First we make the substitutions $v:= u^m$ and $\tilde{g}:= g^m$ to write (\ref{fracPMie}) as 
\begin{equation}\label{vfracPMie}
\left\{
\begin{aligned}
\partial_t(v^\frac{1}{m})+ L^s_{\gamma}v+ \lambda(x)v^\frac{1}{m} &= 0,\quad \Omega\times (0, T),\\
v&= \tilde{g},\quad \Omega_e\times (0, T),\\
v&= 0,\quad \Omega\times \{0\}.\\
\end{aligned}
\right.
\end{equation}
We define the associated Dirichlet-to-Neumann map
\begin{equation}\label{vglDN}
\tilde{\Lambda}_{\gamma,\lambda}: \tilde{g}\to L^s_{\gamma}v|_{\Omega_e\times (0, T)}.
\end{equation}
Clearly, the knowledge of 
$$\tilde{\Lambda}_{\gamma,\lambda}\tilde{g}|_{W_2\times (0, T)},\quad \tilde{g}\in C^\infty_c(W_1)$$
is equivalent to the knowledge of
$$\Lambda_{\gamma,\lambda}g|_{W_2\times (0, T)},\quad g\in S_{W_1}.$$

In the rest of this paper, we use $v^{(h)}$ to denote the solution of (\ref{vfracPMie}) corresponding to
$\tilde{g}:= hg_0(x)$, $g_0\in C^\infty_c(W_1)$ where $h>0$ is a parameter.
Let $m'$ be the constant s.t. $\frac{1}{m}+ 
\frac{1}{m'}= 1$. Choose a constant $\alpha$ s.t. $\alpha> m'-1$. The constants $C, C', C_1,\cdots$ in this section may depend on some parameters but will never depend on $T, h$ and $g_0$.

\subsection{Time-integral transform}
We consider the time-integral transform
\begin{equation}\label{tit}
V(x):= \int^T_0(T-t)^\alpha v(x, t)\,dt.
\end{equation}
This transform was used in \cite{carsteaU2021inverse} to relate the nonlinear local parabolic problem to the linear local elliptic problem. We will see that it still works in our fractional setting.

We also define
$$M(x):= 
\alpha\int^T_0(T-t)^{\alpha-1}v^\frac{1}{m}(x,t)\,dt,\qquad
N(x):= \lambda(x)\int^T_0(T-t)^{\alpha}v^\frac{1}{m}(x,t)\,dt.$$
By H{\"o}lder's inequality we get the pointwise estimate
\begin{equation}\label{ptM}
|M(x)|= |\alpha\int^T_0(T-t)^\frac{\alpha-m'}{m'}((T-t)^\alpha v)^\frac{1}{m}(x,t)\,dt|\leq C_{\alpha,m}T^{\frac{\alpha}{m'}-\frac{1}{m}}|V(x)|^\frac{1}{m}.
\end{equation}
Similarly we can show that
\begin{equation}\label{ptN}
|N(x)|\leq C'_{\alpha,m}T^{\alpha+m'}|\lambda(x)||V(x)|^\frac{1}{m}.
\end{equation}

Note that by Proposition 3.3 for each $t\in (0,T)$, $$\tilde{w}^{(h)}(t):=
v^{(h)}(t)- hg_0\in \tilde{H}^s(\Omega),\qquad (\tilde{w}^{(h)})^\frac{1}{m},
\partial_t((\tilde{w}^{(h)})^\frac{1}{m})\in L^2(0, T; H^{-s}(\Omega))$$
so by the equation in (\ref{vfracPMie}) we get 
$$L^s_{\gamma}\tilde{w}^{(h)}|_{\Omega\times (0, T)}\in L^2(0, T; H^{-s}(\Omega))$$
and thus we get
$$\tilde{w}^{(h)}\in L^2(0, T; \tilde{H}^s(\Omega)),\qquad V^{(h)}\in H^s(\mathbb{R}^n)$$
so $L^s_{\gamma}V^{(h)}\in H^{-s}(\mathbb{R}^n)$.

Moreover, by applying (\ref{tit}) to (\ref{vfracPMie}) and integration by parts with respect to $t$ we get
\begin{equation}\label{ellMN}
\left\{
\begin{aligned}
L^s_{\gamma}V^{(h)} &= M^{(h)}+ N^{(h)},\quad x\in\Omega,\\
V^{(h)}&= C_\alpha T^{1+\alpha}hg_0,\quad\,\, x\in\Omega_e.\\
\end{aligned}
\right.
\end{equation}
By (\ref{ptM}) and (\ref{ptN}) we get the $L^2$-estimates
\begin{equation}\label{L2MN}
||M^{(h)}||_{L^2(\Omega)}\leq C''_{\alpha,m}T^{\frac{\alpha}{m'}-\frac{1}{m}}
||V^{(h)}||_{L^2(\Omega)}^\frac{1}{m},\quad 
||N^{(h)}||_{L^2(\Omega)}\leq C''_{\alpha,m,\lambda}T^{\alpha+ m'}
||V^{(h)}||_{L^2(\Omega)}^\frac{1}{m}.
\end{equation}
By (\ref{ellMN}) and Proposition 2.1 we get the estimate
\begin{equation}\label{ellHs}
||V^{(h)}||_{H^s(\mathbb{R}^n)}\leq C(||M^{(h)}+ N^{(h)}||_{H^{-s}(\Omega)}
+ C_\alpha T^{1+\alpha}h||g_0||_{H^s(\mathbb{R}^n)}).
\end{equation}
Combining (\ref{ellHs}) with (\ref{L2MN}), we get
$$||V^{(h)}||_{H^s}\leq C_1
((T^{\frac{\alpha}{m'}-\frac{1}{m}}+ T^{\alpha+m'})||V^{(h)}||_{H^s}^\frac{1}{m}+ T^{1+\alpha}h||g_0||_{H^s}).$$

We can assume 
$$1- C_1(T^{\frac{\alpha}{m'}-\frac{1}{m}}+ T^{\alpha+m'})\geq \frac{1}{2}.$$ 
(Otherwise we just replace $T$ by a smaller $T'$ in (\ref{vfracPMie}).) Then we get
$$||V^{(h)}||_{H^s}\leq \max\{1,\, 2C_1T^{1+\alpha}h||g_0||_{H^s}\}$$
so for $g_0\neq 0$ and sufficiently large $h$ (depending on the norm of $g_0$) we get
\begin{equation}\label{HsV}
||V^{(h)}||_{H^s}\leq 2C_1T^{1+\alpha}h||g_0||_{H^s}.
\end{equation}

Now we consider the equality
\begin{equation}\label{V0Rid}
V^{(h)}= C_\alpha T^{1+\alpha}hV_0+ R^{(h)}
\end{equation}
where $V_0$ is the solution of the problem
\begin{equation}\label{ellV0}
\left\{
\begin{aligned}
L^s_{\gamma}V_0 &= 0,\quad\,\,\, x\in\Omega,\\
V_0&= g_0,\quad x\in\Omega_e\\
\end{aligned}
\right.
\end{equation}
and $R^{(h)}$ is the solution of the problem 
\begin{equation}\label{ellRh}
\left\{
\begin{aligned}
L^s_{\gamma}R^{(h)} &= M^{(h)}+ N^{(h)},\quad x\in\Omega,\\
R^{(h)}&= 0,\qquad\qquad\qquad x\in\Omega_e.\\
\end{aligned}
\right.
\end{equation}

We apply $L^s_{\gamma}$ to (\ref{V0Rid}) to obtain
\begin{equation}\label{hV0Rid}
h^{-1}L^s_{\gamma}V^{(h)}= C_\alpha T^{1+\alpha}L^s_{\gamma}V_0+ h^{-1}L^s_{\gamma}R^{(h)}.
\end{equation}
Note that (\ref{ellRh}) and (\ref{L2MN}) imply
$$||R^{(h)}||_{H^s}\leq C||M^{(h)}+ N^{(h)}||_{H^{-s}(\Omega)}
\leq C_1(T^{\frac{\alpha}{m'}-\frac{1}{m}}+ T^{\alpha+m'})||V^{(h)}||_{H^s}^\frac{1}{m}.$$
Hence for each $g_0\neq 0$, by (\ref{HsV}) we get
\begin{equation}\label{Oh}
h^{-1}L^s_{\gamma}R^{(h)}= O(h^{\frac{1}{m}-1})
\end{equation}
as $h\to \infty$ in $H^{-s}$-norm. 

Now the asymptotic behavior
of $h^{-1}L^s_{\gamma}V^{(h)}$ as $h\to \infty$ is clear from $(\ref{hV0Rid})$.

\subsection{Proof of the main theorem}
Before we prove the main theorem, we first need to present 
the uniqueness result for the linear fractional elliptic problem.
Recall that the associated exterior Dirichlet problem is given by (\ref{eDp}) and the Dirichlet-to-Neumann map is defined by (\ref{DNlin}). 

The following proposition is a combination of Theorem 1.4 and Theorem 1.5 in \cite{ghosh2021calder}, which is a fractional analogue of the uniqueness result for the classical Calder\'on problem.

\begin{prop}
Let $W_1, W_2\subset \Omega_e$ be nonempty and open.
Suppose
$0< \gamma^{(1)}, \gamma^{(2)}\in C^\infty(\mathbb{R}^n)$ and $\gamma^{(1)}= \gamma^{(2)}= 1$ in $\Omega_e$.
Suppose 
\begin{equation}
\Lambda^{lin}_{\gamma^{(1)}}g|_{W_2}
= \Lambda^{lin}_{\gamma^{(2)}}g|_{W_2}
\end{equation}
for all $g\in C^\infty_c(W_1)$. Then $\gamma^{(1)}= \gamma^{(2)}$
in $\Omega$.
\end{prop}

We are ready to prove Theorem 1.1. We will first determine $\gamma$ based on the asymptotic analysis in the previous subsection.

\begin{proof}
(\textbf{Uniqueness of $\gamma$}) By the assumption (\ref{IdDN}) we get 
$$\tilde{\Lambda}_{\gamma^{(1)},\lambda^{(1)}}\tilde{g}|_{W_2\times (0, T)}
= \tilde{\Lambda}_{\gamma^{(2)},\lambda^{(2)}}\tilde{g}|_{W_2\times (0, T)},\qquad \tilde{g}\in C^\infty_c(W_1),$$
which implies
$$ L^s_{\gamma^{(1)}}v^{(h)}_1|_{W_2\times (0, T)}
= L^s_{\gamma^{(2)}}v^{(h)}_2|_{W_2\times (0, T)}.$$
Now we apply the time-integral transform (\ref{tit})
to the equality above to obtain
$$ L^s_{\gamma^{(1)}}V^{(h)}_1|_{W_2}= L^s_{\gamma^{(2)}}V^{(h)}_2|_{W_2}.$$
Let $h\to \infty$. By (\ref{hV0Rid}) and (\ref{Oh}) we get
$$L^s_{\gamma^{(1)}}(V_0)_1|_{W_2}= L^s_{\gamma^{(2)}}(V_0)_2|_{W_2},$$
$$i.e.\qquad\Lambda^{lin}_{\gamma^{(1)}}g_0|_{W_2}= \Lambda^{lin}_{\gamma^{(2)}}g_0|_{W_2},\quad g_0\in C^\infty_c(W_1).$$
Now we conclude that $\gamma^{(1)}= \gamma^{(2)}=: \gamma$ in $\Omega$ by Proposition 4.1.
\end{proof}

In \cite{carsteaU2021inverse}, the authors used a more careful asymptotic analysis to determine $\lambda$ after they determined $\gamma$. In our case, the determination of $\lambda$ becomes much simpler due to the unique continuation property of the fractional operator. 

\begin{proof}
(\textbf{Uniqueness of $\lambda$}) We pick a nonzero $g$ in (\ref{fracPMie}).
Now for each $t\in (0, T)$, by (\ref{IdDN}) we get
$$L^s_{\gamma}u^m_1(t)= L^s_{\gamma}u^m_2(t),\qquad x\in W_2.$$
Also note that $$u_1(t)= u_2(t)= g$$ in $\Omega_e$ so we get
$$L^s_{\gamma}(u^m_1(t)- u^m_2(t))= u^m_1(t)- u^m_2(t)= 0$$
in $W_2$, which implies $u_1= u_2:= u$ in $\mathbb{R}^n\times (0, T)$
by Proposition 2.2.

We fix a point $x_0\in \Omega$. We claim that there does not exist an open $U\subset \Omega$ containing $x_0$ s.t. $u= 0$ in $U\times (0, T)$. In fact, suppose such $U$ exists. Then $\partial_t u= 0$ in $U\times (0, T)$
and thus $$L^s_{\gamma}(u^m)= u^m= 0$$ 
in $U\times (0, T)$. However, Proposition 2.2 implies that $u= 0$ in $\mathbb{R}^n\times (0, T)$, which contradicts $g\neq 0$. 

Hence we can choose a sequence $\{x_k, t_k\}$ s.t. 
$$x_k\to x_0,\qquad (x_k, t_k)\in \Omega\times (0,T)$$ and $u(x_k, t_k)\neq 0$ so
$$\lambda^{(j)}(x_0)= \lim_k \lambda^{(j)}(x_k)
= -\lim_k\frac{\partial_tu(x_k, t_k)+ L^s_{\gamma}(u^m)|_{(x_k, t_k)}}{u(x_k, t_k)},\qquad j= 1,2.$$
We conclude that $\lambda^{(1)}= \lambda^{(2)}$ in $\Omega$ since $x_0$ is arbitrary.

\end{proof}
We remark that this method was first used to determine the potential in the fractional Calder\'on problem via a single measurement. See \cite{ghosh2020uniqueness} for more details.

\bibliographystyle{plain}
{\small\bibliography{Reference6}}
\end{document}